\documentclass[12pt]{article}
\usepackage{makeidx}
\usepackage{amssymb}
 \usepackage{latexsym}
\usepackage{tikz}
\usepackage{pgfplots}
\usepackage{pgfkeys}%
\usetikzlibrary{plotmarks}
\usetikzlibrary{calc}
 \usepackage{array}
\usepackage{arabtex}
 
 \usepackage{tabularx}
\usepackage{epsfig}
\usepackage{array}
\usepackage[ansinew]{inputenc}
\usepackage{amsfonts,amsmath}
\usepackage{latexsym}
\usepackage{graphics,color,graphicx,shortvrb}
 
 \usepackage{tabularx}
\newtheorem{Theorem}{Theorem}[section]
\newtheorem{Definition}[Theorem]{Definition}
\newtheorem{Proposition}[Theorem]{Proposition}

\definecolor{ProcessBlue}{cmyk}{1,0,0,0.25}
\definecolor{Black}{cmyk}{0,0,0,1}
\definecolor{Red}{cmyk}{0,1,1,0}
\definecolor{Green}{cmyk}{0.9,0,1,0}
\definecolor{Orange}{cmyk}{0,0.61,0.87,0.1}
\definecolor{Fuchsia}{cmyk}{0.47,0.91,0,0.06}
\definecolor{PineGreen}{cmyk}{0.92,0,0.59,0.25}

\usepackage[plainpages=false,hyperindex=false]{hyperref}
\hypersetup{pdfpagemode=FullScreen}

\begin{document}
 \sloppy
 
\begin{center} 
{\Large \bf  Cournot Maps for Intercepting Evader Evolutions by a 
Pursuer} 
\\  \mbox{}\\ \textbf{\emph{Jean-Pierre Aubin}}\footnote{VIMADES (Viabilité, 
Marchés, Automatique, 
Décisions), 14, rue Domat, 75005, Paris, France \\
aubin.jp@gmail.com,  
http://vimades.com/aubin/},\footnote{\textbf{Acknowledgments} 
\emph{This work was partially supported by the Commission of the 
European Communities under the 7th Framework Programme Marie 
Curie Initial Training Network (FP7-PEOPLE-2010-ITN),  project 
SADCO, contract number 264735 and ANR-11-ASTRID-0041-04.}} and  
\textbf{\emph{Chen Luxi}}\footnote{Université Panthéon-Sorbonne 
and Société VIMADES (Viabilité, Marchés, Automatique et 
Décision), e-mail: clxshd@gmail.com} 
\\   \mbox{} January 6, 2014
\end{center}


\subsubsection*{Abstract}
\emph{Instead of studying evolutions governed by an evolutionary 
system starting at a \emph{given} initial state  on a 
\emph{prescribed future}  time interval, finite or infinite, we 
tackle the problem of looking \emph{both} for a \emph{past 
interval} $[T-D,T]$ of aperture (or length, duration) $D$ and for 
the \emph{viable evolutions arriving at a prescribed terminal 
state} at the end of the temporal window (and thus telescoping if 
more than one such evolutions exist). \\ Hence, given time and 
duration dependent  evolutionary system and viability 
constraints, as well as time dependent departure constraints, the 
Cournot map  associates with any \emph{terminal time $T$ and 
state $x$}  the apertures $D(T,x)$ of the intervals 
$[T-D(T,x),T]$, the starting (or initial) states at the beginning 
of the temporal window from which at least one viable evolution 
will reach the given terminal state $x$ at $T$.  Cournot maps can 
be used by a pursuer to intercept an evader's evolution in 
dynamic game theory. After providing some properties of Cournot 
maps are next investigated, above all, the \emph{regulation map} 
piloting the viables evolutions at each time and for each 
duration from the beginning of the temporal window up to terminal 
time. \\The next question investigated is the selection of 
controls or regulons in the regulation map  whenever several of 
them exist. Selection processes are either \emph{time dependent}, 
when the selection operates at each time, duration and state for 
selecting a regulon satisfying required properties (for instance, 
minimal norm, minimal speed), or \emph{intertemporal}. In this 
case, viable evolutions are required to  optimize some prescribed 
intertemporal functional, as in optimal control. This generates 
value functions, the topics of the second part of this study.}

\subsubsection*{Mathematics Subject Classification:} 
34A60, 90B10, 90B20, 90B99, 93C10, 93C30, 93C99 
 
\subsubsection*{Keywords:} 
Cournot uncertainty, Evader interception,  duration 
structuration, viability,   terminal conditions, intertemporal 
optimality

\section{Introduction}

We attempt to translate mathematically  an important concept of 
\emph{uncertainty} suggested in \emph{Exposition de la théorie 
des chances et des probabilités}, \cite[Cournot]{Cournot}, 1843, 
by \glossary{Cournot (Augustin) [1811-1877]} \emph{Augustin 
Cournot} as the meeting of two independent causal series: 
\emph{``A myriad partial series can coexist in time: they can 
meet, so that a single event, to the production of which several 
events took part, come from several distinct series of generating 
causes.''}  The search for causes amounts in this case to look 
for ``retrodictions'' (so to speak) instead of 
predictions\footnote{This idea probably goes back to the 
presocratic Greeks, according to the biologist \glossary{Danchin 
(Antoine) [1944-]}  \emph{Antoine Danchin} in  
\cite[Danchin]{Danchin3,Danchin4,Danchin5} and his book, \emph{La 
Barque de Delphes. Ce que révèle le texte des génomes}, 
\cite[Danchin]{Danchin2}. He denotes what we called ``Cournot 
uncertainty'' as ``contingent uncertainty'', whereas we use the 
adjective ``contingent'' in viability theory for translating 
mathematically  the uncertainty encapsulated in differential 
inclusions, for actually capturing the concept of redundancy, not 
only describing the telescoping of evolutions, but the choice of 
adequate (for instance, viable), regulons or controls in a 
``contingent reservoir'', which may itself evolve. }.

We suggest to combine this Cournot approach uncertainty with the 
Darwinian view of \emph{contingent uncertainty} (differential 
inclusions) for facing \emph{necessity}\footnote{In 
\cite[Danchin]{Danchin4}, the author quotes   the following 
sentence of \glossary{Leucippus [460-370 av. J.-C.]} 
\emph{Leucippus}: \emph{``Nothing happens in vain, but everything 
from reason (logos), and by necessity''}: ``law'' is described by 
a differential inclusion, ``necessity'' by constraints to abide 
to. This ``Law and Necessity'' statement  is more in tune with 
viability theory  than  \emph{Chance and Necessity},  title  of 
the celebrated book \cite[Monod]{Monod}  by \glossary{Monod 
(Jacques) [1910-1976]} \emph{Jacques Monod}, who attributed this 
concept to \glossary{Democritus [460-370 av. J.-C.]} 
\emph{Democritus}. Indeed, ``chance'' remains to be defined (see 
Chapter~2 of \emph{La valeur n'existe pas. À moins que ...}, 
\cite[Aubin]{valeur} and Chapter~8 of \emph{La  mort  du  devin, 
l'émergence du  démiurge. Essai sur la contingence, la viabilité 
et l'inertie des systèmes}, \cite[Aubin]{mded}). } (viability 
constraints) by introducing the concept of Cournot map.   

We   provide a viability characterization of   Cournot maps which 
relates them to the concept of capture basins viable in an 
environment  (see Chapter~8, p. 273, of \emph{Viability Theory.  
New Directions}, \cite[Aubin, Bayen \&  
Saint-Pierre]{absp}\footnote{In this book, the evolutions are 
still defined on the usual future interval  $[0,T]$  with 
prescribed finite or infinite horizon.}).

More generally, Cournot maps are  motivated by traffic congestion 
(where the duration is the travel time), by economic dynamics 
(where the duration of investment evolution), by population 
dynamics\footnote{See, among an abundant literature, 
\cite[Anita]{Anita}, \cite[Aubin]{RegulationBirths-11},   
\cite[Iannelli]{im95aspd},   \cite[Keyfitz N. \& Keyfitz 
B.]{ke297dem}, \cite[Von Foerster]{foe59dem}, 
\cite[Webb]{web85age}, etc.} (where the duration is age), by 
collision problems  (where the duration is time until collision): 
Cournot maps can also be used by a pursuer to intercept an 
evader's evolution in dynamic game theory.  

Following the suggestion to study evolutions on (sliding) 
\emph{temporal window} $[T-\Omega ,T]$ on which the evolution is 
defined\footnote{see Chapter~5, p. 67, of \emph{Time and Money. 
How Long and How Much Money is Needed to Regulate  a Viable  
Economy},\cite[Aubin]{TM} and \emph{La valeur n'existe pas. À 
moins que ...}, \cite[Aubin]{valeur}.}, \emph{the search for the 
temporal window  is \emph{also} part of the solution of the 
problem}. 

\section{Cournot Maps}

\subsection{Definitions}

Let us consider 

\begin{enumerate}   \item    

a set-valued map $F: \mathbb{R}^{} \times \mathbb{R}^{} \times  
X  \leadsto X$ with which we associate 
\begin{enumerate}   

\item   the arrival map 
$$\mathcal{A}_{F}: \mathbb{R}^{} \times \mathbb{R}^{}_{+}  \times X  
\times X  \leadsto \mathcal{C}(  -\infty,+\infty;X) 
$$ 
associating with any terminal pair $(T,x)$,  aperture $\Omega 
\geq 0$ and $s $ the (possibly empty) set 
$\mathcal{A}_{F}(T,x)[\Omega,s]$ of evolutions $x(\cdot)$ 
restricted to the temporal window $[T-\Omega ,T]$ governed by the 
\index{duration-structured differential inclusion} 
\emph{duration-structured differential inclusion} 

\begin{equation} \label{e:LX-DifIncl}
\forall \;  t \in [T-\Omega ,T], \; \;  x'(t) \; \in \;  
F(t,t-(T-\Omega),x(t)) 
\end{equation}
\emph{defined on the temporal window} $[T-\Omega ,T]$ starting 
from $s  $ at the beginning of the temporal window and arriving 
at  $x(T) = x \in K(T)$ at its end\footnote{This is a 
\index{retrodiction system} \emph{``retrodiction'' evolutionary 
system} in the sense that, for any evolving present time $T$, we 
study the evolution on the past (or historical) temporal window 
$[T-\Omega ,T]$.}. Such an evolution linking $s$ at time 
$T-\Omega $ to $x$ at time $T$ is called a \index{Cournot 
evolution} \emph{Cournot evolution} at $(T,x)$;

\item   the map $(T,x, \Omega )  \leadsto 
\mathcal{A}_{F}(T,x)[\Omega ] := \bigcup_{s \in C(T-\Omega 
)}\mathcal{A}_{F}(T,x)[\Omega,s]$, the set of evolutions defined 
on the temporal window $[T-\Omega ,T]$ arriving at $x$ at time 
$T$; 

\item   the map $(T,x) \leadsto  \bigcup_{ \Omega \geq 0 
}\mathcal{A}_{F}(T,x)[\Omega]$, the set of evolutions   arriving 
at $x$ at time $T$; 
\end{enumerate} 

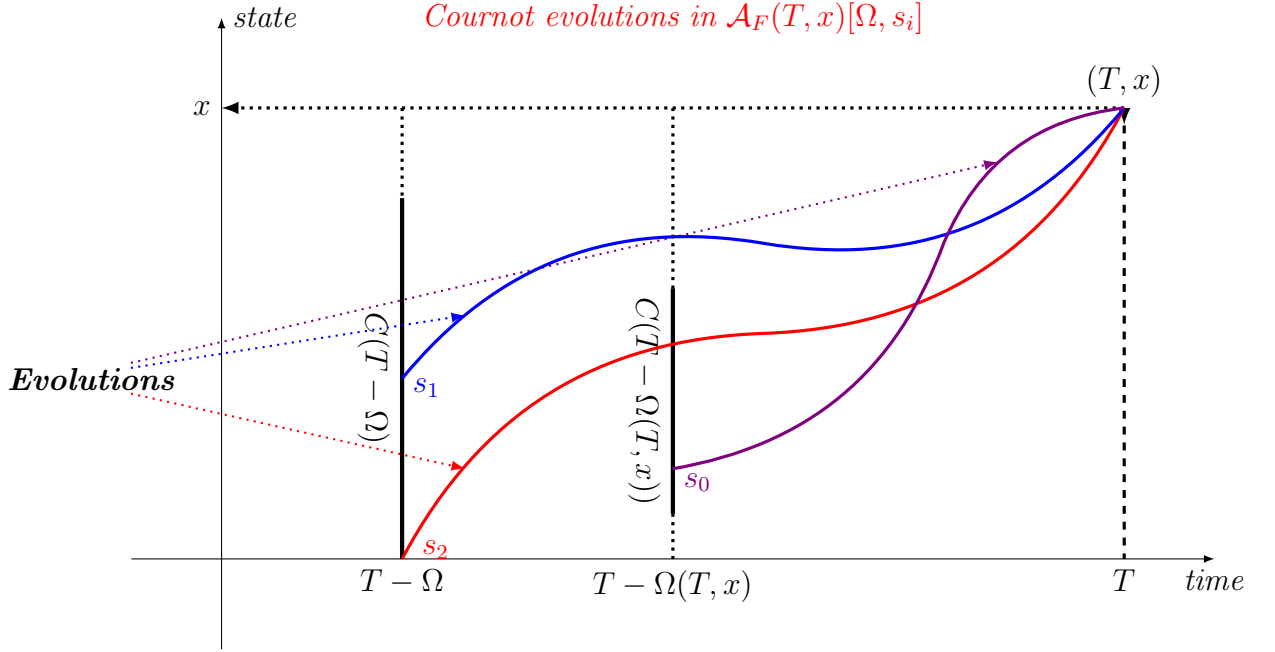
\begin{center}
\begin{tikzpicture}  [scale = 1.2]
\node[text width=7cm, red, text centered]  (Z) at (5,6 ) 
{\emph{Cournot evolutions in \textbf{$\mathcal{A}_{F}(T,x)[\Omega 
,s_{i}]$}}};  
 
\draw (10,5) node[above]{$(T,x)$}; 

\draw[->,>=latex, ] (-1,0) -- (11,0);

\draw[->,>=latex, ] (0, -1) -- (0,6);

\draw (10,0) node[below]{$T$};

\draw (0,5) node[left]{$x$}; 

\draw (2,0) node[below]{$T-\Omega $};

\draw (5,0) node[below]{$T-\Omega (T,x)$};

\draw (11,0) node[below]{\emph{time}};

\draw (0,6) node[right]{\emph{state}}; 

\draw[dotted, very thick] (2,5) -- (2,0); 

\draw[ultra thick] (2,4) -- (2,0) node[midway,below,sloped] {$ 
C(T-\Omega ) $}; 

\draw[dotted, very thick] (5,5) -- (5,0); 

\draw[ultra thick] (5,3) -- (5,.5) node[midway,below,sloped] {$ 
C(T-\Omega(T,x) ) $};     

\draw[->,>=latex, dotted, very thick] (10,5) -- (0,5); 

\draw[>-,>=latex, dashed, very thick] (10,5) -- (10,0);

\node[text width=1.5cm, text centered]  (B) at (-1.75,2) 
{\emph{\textbf{Evolutions}}}; 

\draw[->,>=latex,dotted, thick,  red] (B) -- (2.7,1); 

\draw[->,>=latex,dotted, thick, blue] (B) -- (2.7,2.7); 

\draw[->,>=latex,dotted, thick,  violet] (B) -- (8.6,4.4); 


\draw[very thick, red]   (2,0)  to[bend left]  (6,2.5) to[bend 
right]  (10,5);

\draw (2.1,0.1) node[red, right]{$s_{2}$};

\draw[very thick, blue]   (2,2)  to[bend left]  (6,3.5) to[bend 
right]  (10,5); 

\draw (2,1.9) node[blue, right]{$s_{1}$};

\draw[very thick, violet]   (5,1)  to[bend right]  (8,3.5) 
to[bend left]  (10,5); 

\draw (5,0.85) node[violet, right]{$s_{0}$};

\end{tikzpicture} \end{center} 

\item a (duration dependent) \index{environmental tube} 
\emph{environmental tube}  $\mathbb{K} : (t,d) \in  \mathbb{R}^{} 
\times \mathbb{R}^{}_{+} \leadsto K(t,d) \subset  X$  in which 
evolutions are required to be viable;

\item   a \index{departure tube} \emph{departure 
tube}\footnote{The departure sets $C(d)$ can be empty for some 
departure dates $d$. If, for instance, the beginning $d_{fix}$ of 
the temporal window is prescribed and not computed, then the 
departure sets $C(d)=\emptyset$ are empty for all $d \ne 
d_{fix}$. If we want that all departure dates are later than a 
date $d_{min}$, we assume that the departure sets 
$C(d)=\emptyset$ are empty for all $d < d_{min}$. The use of 
departure maps cover many different situations. In ethology, 
departure maps could translate mathematically the 
\glossary{Lorenz (Konrad) [1903-1989]} \emph{Konrad Lorenz} 
imprinting (actually discovered by the 19th-century  biologist 
\glossary{Spalding (Douglas Alexander) [1841-1877]} \emph{Douglas 
Spalding}, rediscovered by  \glossary{Heinroth (Oskar) 
[1871-1945]} \emph{Oskar Heinroth}, Lorenz' mentor), associating 
with given dates of cognitive development the perception of the 
environment triggering imprinted behaviors (such as recognition 
of the mother, etc.).}  $\mathbb{C}: t\in  \mathbb{R}^{}  
\leadsto  C(t) \subset  K(t,0)$, which associates with any 
(departure)  date $t\in \mathbb{R}^{}$ the subset  of   state $s 
\in C(t)$ from which evolutions start at time $t$.

\end{enumerate}

 The question arises whether we can find the subset of those 
initial states.

\begin{Definition} 
\symbol{91}\textbf{Cournot Map}\symbol{93}\label{}\index{} The 
\index{Cournot map} \emph{Cournot map} $\mbox{\rm 
Cour}_{F}(\mathbb{K},\mathbb{C}): \mbox{\rm Graph}(K)  \leadsto 
\mathbb{R}^{}_{+} \times X $ of the departure tube $\mathbb{C}$ 
viable in the environmental tube $\mathbb{K}$ under the 
differential inclusion (\ref{e:LX-DifIncl}), p. 
\pageref{e:LX-DifIncl}, is the set-valued map   associating  with 
any terminal pair $(T,x) \in \mbox{\rm Graph}(\mathbb{K})$ the 
subset  $\mbox{\rm Cour}_{F}(\mathbb{K},\mathbb{C})(T,x)$ of 
pairs $(\Omega ,s) \in \mathbb{R}^{}_{+} \times X$ such that 
\begin{enumerate}   
\item  $s \in C(T-\Omega )$;  

\item   there exists at least one Cournot evolution $x(\cdot) \in 
\mathcal{A}_{F}(T,x)[\Omega,s]$ starting from $x(T-\Omega )=s$ at 
time $T-\Omega $, arriving at $x(T)=x$ at time $T$ and  viable in 
the environmental tube $x(\cdot)$ on the temporal window 
$[T-\Omega ,T]$ of aperture $\Omega \geq 0$ in the sense that 

\begin{equation} \label{e:}   
\forall \;  t \in [T-\Omega ,T], \; \; x(t) \; \in \; K(t, 
t-(T-\Omega))
\end{equation}
\end{enumerate}

When there is no environmental constraint, we simply set 
$\mbox{\rm Cour}_{F}(\mathbb{C})$. \\  We stress the fact that we 
look for  \emph{both a temporal window} $[T-\Omega ,T]$ 
\textbf{and} a \emph{viable evolution} $x(\cdot) $ governed by 
$x'(t)\in F(t,t-(T-\Omega),x(t)) $ on this temporal window.  
\end{Definition} 

The concept of Cournot map encapsulates several features. The 
first one is the concept of Cournot (or minimal aperture):

\begin{Definition} 
\symbol{91}\textbf{Cournot Aperture}\symbol{93}\label{}\index{} 
The \index{Cournot aperture} \emph{Cournot aperture} function
$\Omega_{F}(\mathbb{K},\mathbb{C}) $ associates with any  time 
$T$ and at arrival state $x$ the smallest aperture 
\begin{equation} \label{e:}  
\Omega_{F}(\mathbb{K},\mathbb{C})(T,x) \; := \; \inf_{(\Omega,s) 
\; \in \; \mbox{\rm 
Cour}_{F}(\mathbb{K},\mathbb{C})(T-\Omega_{F}(\mathbb{K},\mathbb{C})(T,x),x) 
}\Omega 
\end{equation} 
of the Cournot temporal window 
$[T-\Omega_{F}(\mathbb{K},\mathbb{C})(T,x),T]$. Its inverse 
$\displaystyle{\frac{1}{\Omega_{F}(\mathbb{K},\mathbb{C})(T,x) 
}}$ is called the \index{viability liquidity} \emph{Cournot 
liquidity} in economics.\\ When there is ambiguity, we set 
$\Omega (T,x) :=\Omega_{F}(\mathbb{K},\mathbb{C}) (T,x)$
\end{Definition}

Next, we extract from the knowledge of the Cournot map the 
arrival tube of arrival dates and states which can be reached:

\begin{Definition} 
\symbol{91}\textbf{Cournot Arrival Tubes}\symbol{93} 
\label{d:ArrivalTube}\index{} The Cournot map $\mbox{\rm 
Cour}_{F}(\mathbb{K},\mathbb{C})$ generates 
 the \index{arrival tube} \emph{arrival tube} $\mbox{\rm 
Aval}_{F}(\mathbb{K},\mathbb{C}) \subset \mathbb{K}$ defined by  
\begin{equation} \label{e:ArrivalTube}   
\mbox{\rm Graph}( \mbox{\rm Aval}_{F}(\mathbb{K},\mathbb{C}) ) \; 
:= \; \mbox{\rm Dom}(\mbox{\rm Cour}_{F}(\mathbb{K},\mathbb{C}))
\end{equation}
which associates with any $ T \in \mathbb{R}^{}$ the (possibly 
empty) subset $  \mbox{\rm Aval}_{F}(\mathbb{K},\mathbb{C})(T)$  
of arrival states $x \in K(T)$ at which at least one 
\emph{viable} Cournot evolution starting from some $s \in 
C(T-\Omega )$ at $T-\Omega $ for some aperture $\Omega \geq 0$ 
arrives at $x$ at time $T$.   
\end{Definition}

In other words, the set-valued map $F$ generates the map 
\begin{equation} \label{e:}   
\mathbb{C}  \mapsto  \mbox{\rm Aval}_{F}(\mathbb{K},\mathbb{C}) 
 \;\subset\;  \mathbb{K} 
\end{equation}
\emph{mapping departure tubes $\mathbb{C}$ to arrival tubes 
$\mbox{\rm Aval}_{F}(\mathbb{K},\mathbb{C})$}.

The Cournot map generates in turn the Cournot starting tube 
contained in the departure tube: 

\begin{Definition} 
\symbol{91}\textbf{Cournot Starting Tubes}\symbol{93} 
\label{d:ArrivalTube}\index{} The Cournot map $\mbox{\rm 
Cour}_{F}(\mathbb{K},\mathbb{C})$ generates the \index{starting 
map} \emph{starting map} $\mbox{\rm 
Start}_{F}(\mathbb{K},\mathbb{C}) \subset \mathbb{C}$ defined by 
\begin{equation} \label{e:}   
(T,x)  \leadsto \mbox{\rm Start}_{F}(\mathbb{K},\mathbb{C})(T,x) 
\; := \; \left\{ (T-\Omega ,s) \right\}_{(\Omega,s) \in \mbox{\rm 
Cour}_{F}(\mathbb{K},\mathbb{C})(T,x)}   
\end{equation}
We denote by $\mbox{\rm Start}_{F}(\mathbb{K},\mathbb{C})(T,x) 
[\Omega ]\; := \; \mbox{\rm Cour}_{F}(\mathbb{C})(T,x)[\Omega]$ 
the subset defined by 
\begin{equation} \label{e:} \mbox{\rm 
Start}_{F}(\mathbb{K},\mathbb{C})(T,x) [\Omega ]  \; := \;  
\left\{ s \; \mbox{ such that} \; (\Omega ,s) \; \in \; \mbox{\rm 
Cour}_{F}(\mathbb{C})(T,x)\right\}
\end{equation}
providing the starting states $s \in C(T-\Omega )$. In 
particular, we single out the  \index{Cournot earliest starting 
map} \emph{Cournot earliest starting map}
\begin{equation} \label{e:}   
(T,x)   \leadsto \mbox{\rm 
Cour}_{F}(\mathbb{C})(T,x)[\Omega_{F}(\mathbb{K},\mathbb{C}(T,x)]  
\;\subset\; C(T-\Omega_{F}(\mathbb{K},\mathbb{C})(T,x) ) 
\end{equation}
associating the set of starting states $s \in  
C(T-\Omega_{F}(\mathbb{K},\mathbb{C})(T,x) )$ from which starts a 
viable Cournot evolution $x(\cdot) \in 
\mathcal{A}_{F}(T,x)[\Omega_{F}(\mathbb{K},\mathbb{C})(T,x),s]$ 
reaching $x$ at time $T$ with minimal duration.\\
The tube $\mbox{\rm Start}_{F}(\mathbb{K},\mathbb{C})(\mbox{\rm 
Aval}_{F}(\mathbb{K},\mathbb{C})) \subset  \mathbb{C}$ is the 
\index{starting tube} \emph{starting tube} of the Cournot map, 
the subset of starting times and states $(d,s)$ from which at 
least one evolution arrives at some $(T,x) \in \mbox{\rm 
Aval}_{F}(\mathbb{K},\mathbb{C})$ in the arrival tube. 
 \end{Definition} 

\subsection{Cournot Tube of a Pursuer for Intercepting an Evader Evolution}

Cournot tubes could be of some use in the context of 
pursuer-evader dynamical games. We consider the problem from the 
point of view of the pursuer,   who has computed its Cournot map 
$\mbox{\rm Cour}_{F}(\mathbb{K},\mathbb{C})$.

Assume that at time $t_{0}$, the pursuer observes an 
evolution\footnote{This evolution can be extrapolated from the 
knowledge of the evolution on an adequate interval $[t_{0}-\Omega 
_{0}, t_{0} ]$ or, knowing the dynamics of the evader, an 
evolution starting at $\xi_{0}(t_{0})$ governed by the evader's 
dynamics.} $\xi_{0}(\cdot) : \symbol{91} t_{0}, +\infty 
\symbol{91} \mapsto X$.  

Can the pursuer intercept the evolution $\xi_{0}(\cdot)$, and, if 
the answer is positive, when and how?  Cournot maps can be used 
for answering these questions\footnote{The literature on 
differential games from \emph{Differential games}, 
\cite[Isaacs]{isaacs65} is so abundant that it is impossible to 
quote all the contributions, which figure, for instance, in the 
recent proceedings, \emph{Advances in Dynamic Games: Theory, 
Applications, and Numerical Methods for Differential and 
Stochastic Games}, \cite[Cardaliaguet  \& 
Cressman]{CardaliaguetCressman}. However, viability techniques 
have been introduced in Chapter~14 of \emph{Viability Theory}, 
\cite[Aubin]{avt}, \cite[Cardaliaguet \& 
Plaskacz]{CardaliaguetPlaskacz}, \cite[Cardaliaguet, Quincampoix  
\& Saint-Pierre]{cqp93avk,cqsp00dg}  among many other articles. }.

\begin{Theorem} 
\symbol{91}\textbf{Capturability of the Evader 
Evolution}\symbol{93}\label{}\index{} Let us assume that the 
arrival tube $ \mbox{\rm Aval}_{F}(\mathbb{K},\mathbb{C})$ of the 
Cournot map of the pursuer is closed. We associate with it and 
with the evader evolution $\xi_{0}(\cdot)$ the capturability 
state $ (T^{\flat}_{\xi_{0}} , x^{\flat}_{\xi_{0}})$   defined by

\begin{equation} \label{e:} \left\{ \begin{array}{ll} 
(i) & T^{\flat}_{\xi_{0}} \; := \;  \inf_{\{t \geq t_{0}  \; 
\mbox{ such that} \; (t,\xi_{0}(t)) \in \mbox{\rm 
Graph}(\mbox{\rm Aval}_{F}(\mathbb{K},\mathbb{C}))\}}t\\  (ii) & 
x^{\flat}_{\xi_{0}} \; := \; \xi_{0}(T^{\flat}_{\xi_{0}} ) 
\end{array} \right. 
\end{equation} 
If $T^{\flat}_{\xi_{0}}<+\infty$ is finite and if the duration 
$\Omega _{0}:= \Omega(T^{\flat}_{\xi_{0}},x^{\flat}_{\xi_{0}}) 
\leq T^{\flat}_{\xi_{0}}-t_{0}$ is smaller or equal to the 
duration $T^{\flat}_{\xi_{0}}-t_{0}$, then the evader evolution 
$\xi_{0}(\cdot)$ is captured by a viable Cournot evolution 
$x_{0}(\cdot) \in \mathcal{A}_{F} 
(T^{\flat}_{\xi_{0}},x^{\flat}_{\xi_{0}}) [\Omega _{0}, s _{0}]$ 
where $s _{0}=x_{0}(T^{\flat}_{\xi_{0}} -\Omega _{0})$.
\end{Theorem}

\textbf{Proof} --- \hspace{ 2 mm} The case when 
$T^{\flat}_{\xi_{0}}=+\infty$ means that the evader evolution is 
not capturable by the pursuer. Otherwise, the pair 
$(T^{\flat}_{\xi_{0}},x^{\flat}_{\xi_{0}})$ belongs to the graph 
$\mbox{\rm Graph}( \mbox{\rm Aval}_{F}(\mathbb{K},\mathbb{C}) )$ 
of the arrival tube of the pursuer. Therefore, there exist one 
Cournot aperture $\Omega _{0}$, one  starting state $s_{0} \in 
C(T^{\flat}_{\xi_{0}}-\Omega _{0})$  and one viable Cournot 
evolution $x_{0}(\cdot)$ linking $s_{0}$ at time 
$T^{\flat}_{\xi_{0}}-\Omega_{0}$ to $x^{\flat}_{\xi_{0}}:= 
\xi_{0}(T^{\flat}_{\xi_{0}})$ at time $T^{\flat}_{\xi_{0}}$, and 
thus intercepting the evader at time $T^{\flat}_{\xi_{0}}$ since 
$t_{0} \leq T^{\flat}_{\xi_{0}}-\Omega_{0}$ by assumption. \hfill 
$\;\; \blacksquare$ \vspace{ 5 mm}

Naturally, the assumption that the observed evolution 
$\xi_{0}(\cdot)$ at $t_{0}$ is known is too strong, since 
predictions are most of the time doomed to fail. Another 
observation may have to be made at a future time $t_{1} \in 
[T^{\flat}_{\xi_{0}}-\Omega_{0},T^{\flat}_{\xi_{0}}]$. It may 
happen that at time $t_{1}$ starts   another evolution 
$\xi_{1}(\cdot)$. 

In this case, at that time $t_{1}$, the state of the pursuer 
evolution $x_{0}(\cdot) $ is no longer viable in the departure 
tube $\mathbb{C}$. This departure tube has to replaced by the 
tube reduced to  $\{x_{0}(\cdot)\}$ from which the evolution a 
possible correction must be made. 

We thus compute the Cournot map $\mbox{\rm 
Cour}_{F}(\mathbb{K},\{x_{0}(\cdot)\})$ so that, for any $t \in  
[T^{\flat}_{\xi_{0}} -\Omega _{0}, T^{\flat}_{\xi_{0}}]$, 
$x_{0}(t) \in \mbox{\rm 
Cour}_{F}(\mathbb{K},\{x_{0}(\cdot)\})(T^{\flat}_{\xi_{0}}, 
x^{\flat}_{\xi_{0}})[T^{\flat}_{\xi_{0}}-t]$.

We next introduce the pair 

\begin{equation} \label{e:} \left\{ \begin{array}{ll} 
(i) & T^{\flat}_{\xi_{1}}  \; := \;  \inf_{\{t \in 
[t_{1},T^{\flat}_{0}]   \; \mbox{ such that} \; (t,\xi_{1}(t)) 
\in \mbox{\rm Graph}(\mbox{\rm 
Aval}_{F}(\mathbb{K},\{x_{0}(\cdot)\}))\}}t\\  (ii) & 
x^{\flat}_{\xi_{1}} \; := \; \xi_{1}(T^{\flat}_{\xi_{1}} ) 
\end{array} \right. 
\end{equation}

The case when $T^{\flat}_{\xi_{1}}=+\infty$ means that the evader 
evolution is not capturable by the pursuer before 
$T^{\flat}_{\xi_{0}}$. Otherwise  $T^{\flat}_{\xi_{1}} \leq 
T^{\flat}_{\xi_{0}}$. We set  $ \Omega _{1}:=\Omega 
(T^{\flat}_{\xi_{1}}, x^{\flat}_{\xi_{1}})$, we take $s_{1}  \in 
\mbox{\rm Cour}_{F}(\mathbb{K},\{x_{0}(\cdot)\} 
(T^{\flat}_{\xi_{1}}, x^{\flat}_{\xi_{1}}))[\Omega_{1}]$ and a 
viable Cournot evolution $x_{1}(\cdot) \in 
\mathcal{A}_{F}(T^{\flat}_{\xi_{1}}, x^{\flat}_{\xi_{1}})[\Omega 
_{1}, s_{1}]$ linking $s_{1}$ at time $T^{\flat}_{\xi_{1}} 
-\Omega _{1}$ to $\xi_{1}(T^{\flat}_{\xi_{1}})$ at time 
$T^{\flat}_{\xi_{1}}$. 

\emph{If $\Omega _{1} \leq T^{\flat}_{\xi_{1}}-t_{1}$, the new 
evader evolution $\xi_{1}(\cdot)$ can be intercepted by a Cournot 
evolution $x_{1}(\cdot) \in \mathcal{A}_{F}(T^{\flat}_{\xi_{1}}, 
x^{\flat}_{\xi_{1}})[\Omega _{1}, s_{1}]$ at time 
$T^{\flat}_{\xi_{1}}$ since   $t_{1} \leq 
T^{\flat}_{\xi_{s}}-\Omega_{1}$.}

We can reiterate this process until interception happens when the 
last prediction $\xi_{j}$ at time $t_{j} \in   
[T^{\flat}_{\xi_{j-1}}-\Omega_{j},T^{\flat}_{\xi_{j-1}}]$ is 
true. 
  
\mbox{}

 
Since Cournot maps can be characterized in terms of viable 
capture basins, they inherit their properties, among them, the 
ability of computing them thanks to the capture basin algorithm. 
They can be used in the field of pursuer-evader dynamical games. 

\subsection{Properties of Cournot Maps}

We observe how the Cournot map $\mbox{\rm 
Cour}_{F}(\mathbb{K},\mathbb{C})(t,x(t))[t-(T-\Omega)] $ evolves 
along a viable Cournot evolution on the temporal window 
$[T-\Omega,T]$ reaching $x$ at time $T$, an obvious consequence 
of the Bilateral Fixed Point of Capture Basins (see 
Theorem~10.2.5, p. 379, of  \emph{Viability Theory.  New 
Directions}, \cite[Aubin, Bayen \&  Saint-Pierre]{absp}):

\begin{Proposition} 
\symbol{91}\textbf{Evolution of Cournot 
Maps}\symbol{93}\label{}\index{} Let us consider $(T,x) \in 
\mbox{\rm Aval}_{F}(\mathbb{K},\mathbb{C})$. \begin{enumerate}   
\item   For any 
 $(\Omega ,s) \in 
\mbox{\rm Cour}_{F}(\mathbb{K},\mathbb{C})(T,x)$  and any viable 
Cournot evolution $x(\cdot) \in  
\mathcal{A}_{(F,\mathbb{K})}(T,x)[\Omega,s]$ linking $s$ to $x$, 
then,
\begin{equation}  \forall \; t 
\; \in \; [T-\Omega ,T], \; \;  (T-\Omega ,s ) \; \in \; 
\mbox{\rm Cour}_{F}(\mathbb{K},\mathbb{C})(t,x(t))[t-(T-\Omega)]
\end{equation}

\item   For any $(\Omega_{1} ,s_{1}) \in \mbox{\rm 
Cour}_{F}(\mathbb{K},\mathbb{C})(T-\Omega ,s)$  and any viable 
Cournot evolution $x_{1}(\cdot) \in  
\mathcal{A}_{(F,\mathbb{K})}(T,x)[\Omega_{1},s_{1}]$ linking 
$s_{1}$ to $s$, then the concatenation $(x_{1}\diamondsuit 
x)(\cdot)$   is a viable Cournot evolution linking $s_{1}$ to 
$x$, and, $\forall \; t \; \in \; [T-(\Omega +\Omega _{1}) ,T],$
\begin{equation}   \; \;  (T-(\Omega +\Omega 
_{1}) ,s_{1} ) \; \in \; \mbox{\rm 
Cour}_{F}(\mathbb{K},\mathbb{C})(t,(x_{1}\diamondsuit 
x)(t))[t-(T-(\Omega +\Omega _{1}))]
\end{equation}
\end{enumerate}
Consequently,

\begin{equation} \label{e:}   
\forall \; t \in [T-\Omega ,T], \; \; t  \leadsto \in \mbox{\rm 
Cour}_{F}(\mathbb{K},\mathbb{C})(t,x(t))  \; \in  \;  \mbox{\rm 
Cour}_{F}(\mathbb{K},\mathbb{C})(T,x) \;\mbox{\rm and is 
increasing} 
\end{equation}
\end{Proposition}

Cournot evolutions $x(\cdot) \in \mathcal{A}_{F}(T,x)[\Omega ,s]
]$ are  not only  in the in the departure tube (after $\Omega 
(T,x)$, at least), so that for we replace it by the Cournot 
evolution $\{x(\cdot)\}$ it self.

\begin{Proposition} 
\symbol{91}\textbf{Viability Property of Cournot 
Evolutions}\symbol{93}\label{}\index{} For any 
 $(\Omega ,s) \in 
\mbox{\rm Cour}_{F}(\mathbb{K},\mathbb{C})(T,x)$  and any viable 
 Cournot evolution $x(\cdot) \in  
\mathcal{A}_{(F,\mathbb{K})}(T,x)[\Omega,s]$ linking $s$ to $x$, 
 
\begin{equation} \label{e:ViabCourEv}   
\forall \;  t \in [T-\Omega ,T], \; \;  x(t)  \; \in \;  
\mbox{\rm Cour}_{F}(\mathbb{K},\{x(\cdot)\})(T,x)[T-t] 
\end{equation} 
\end{Proposition}

Next, we adapt the \index{dilation property} \emph{dilation 
property} of Cournot maps stating in essence that the Cournot map 
of the union of departure tubes is the union of the Cournot maps 
of these departure tubes. This morphism property plays a crucial 
role in computational issues since it allows a parallelization of 
the computation of the Cournot maps.  

We recall that a hypermap $\mathbb{V}$ is  a \emph{dilation} 
\index{dilation} if  $\displaystyle{\mathbb{V}\left(  \bigcup_{i 
\in \mathbb{I}} K_{i} \right) = \bigcup_{i \in 
\mathbb{I}}\mathbb{V}(K_{i})}$ and that any dilation is 
increasing.

\begin{Proposition} 
\symbol{91}\textbf{Morphism Property of Cournot 
Maps}\symbol{93}\label{}\index{} The map  $(F, \mathbb{C})  
\leadsto \mbox{\rm Cour}_{F}(\mathbb{K},\mathbb{C})$ is a 
dilation:

\begin{equation} \label{e:}   
\mbox{\rm Cour}_{\bigcup_{p \in \mathbb{P}}F_{p}}\left( 
\mathbb{K},\bigcup_{i \in \mathbb{I}}\mathbb{C}_{i}\right) \; = 
\; \bigcup_{p \in \mathbb{P}}\bigcup_{i \in \mathbb{I}} \mbox{\rm 
Cour}_{F_{p}}(\mathbb{K},\mathbb{C}_{i})
\end{equation}
and thus, the map $  (F,\mathbb{C})  \leadsto \mbox{\rm 
Cour}_{F}(\mathbb{K},\mathbb{C})$ is increasing.
\end{Proposition}

\subsection{Viability Characterization of Cournot Maps}
 
Our first task is to provide a viability characterization of 
Cournot maps which allows us to transfer the properties of viable 
capture basins to Cournot maps. 


\begin{Theorem} 
\symbol{91}\textbf{Viability Characterization of Cournot 
Maps}\symbol{93}\label{t:VbyChaCournot}\index{} Let us associate 
with the differential inclusion (\ref{e:LX-DifIncl}), p. 
\pageref{e:LX-DifIncl}, the system

\begin{equation} \label{e:LX-DifIncl2} 
\left\{ \begin{array}{ll}  (i)& \overleftarrow{\tau}'(t) \; = \; 
-1   \\ (ii) &   \overleftarrow{\omega}'(t) \; = \; -1 \\ (iii) 
&  \overleftarrow{x}'(t) \; \in \; -F(\overleftarrow{\tau}(t), 
 \overleftarrow{\omega}(t), \overleftarrow{x}(t))\\ (iv) & 
 \overleftarrow{\sigma}'(t) \; = \; 0
\end{array} \right. 
\end{equation}
We introduce the auxiliary environment $\mathcal{K}:=\mbox{\rm 
Graph}(\mathbb{K})\times X$ and the auxiliary target $\mathcal{C} 
\subset \mathcal{K} $ defined by

\begin{equation} \label{e:}   
(t,d,x,s) \; \in \; \mathcal{C} \;\mbox{\rm if and only if} \;  x 
\; \in \; C(t) ,\; d \; = \;  0 \;\mbox{\rm and}\; s = x
\end{equation}
Then  the graph of the Cournot map  $(T,x)  \leadsto  \mbox{\rm 
Cour}_{F}(\mathbb{K},\mathbb{C})(T,x)$ is equal to subset of 
elements $(T,x,\Omega ,s)$ such that $$(T,\Omega ,x,s) \in 
\mbox{\rm Capt}_{(\ref{e:LX-DifIncl2})}(\mathcal{K}, 
\mathcal{C})$$ Therefore, the Cournot map inherits all the 
properties of viable capture basins.
\end{Theorem}

\textbf{Proof} --- \hspace{ 2 mm} To say   $(T,\Omega ,x,s) \in 
\mbox{\rm Capt}_{(\ref{e:LX-DifIncl2})}( \mathcal{K}, 
\mathcal{C})$ belongs to the capture basin   amounts to saying 
that there exist $t^{\star} \geq 0$ and  one evolution 
$(\overleftarrow{\tau}(\cdot), \overleftarrow{\omega}(t), 
\overleftarrow{x}(\cdot), \overleftarrow{\sigma}(t))$ where

\begin{equation} \label{e:}   
\overleftarrow{\tau}(t) \;  = \; T-t , \; 
\overleftarrow{\omega}(t) \;  = \; \Omega -t, \; 
\overleftarrow{x}(t),   \overleftarrow{\sigma}(t) \;  = \; s \;
\end{equation}
governed by differential inclusion (\ref{e:LX-DifIncl2}), p. 
\pageref{e:LX-DifIncl2}, starting at $(T,\Omega,x,s)$ such that

\begin{equation} \label{e:} \left\{ \begin{array}{ll} 
(i) & (\overleftarrow{\tau}(t^{\star}), 
 \overleftarrow{\omega}(t^{\star}), 
\overleftarrow{\xi}(t^{\star}), 
\overleftarrow{\sigma}(t^{\star})) 
\in  \mathcal{C}\\
 &
\mbox{\rm or} \; t^{\star}=\Omega ,  \;\overleftarrow{x}(\Omega) 
\; \in \; C(\overleftarrow{\tau}(\Omega) ) \subset K(\Omega,0)   
\;\mbox{\rm and}\; s \; = \; \overleftarrow{\xi}(\Omega )\\ (ii) 
& \forall \; t \in [0,\Omega ] , \; \; (\overleftarrow{\tau}(t), 
\overleftarrow{\omega}(t),\overleftarrow{\xi}(t),\overleftarrow{\sigma}(t)) 
\in \mathcal{K}\;\mbox{\rm or}\; \overleftarrow{\xi}(t) \in  
K(\overleftarrow{\tau}(t), \overleftarrow{\omega}(t) )
\end{array} \right. 
\end{equation}
 
\mbox{}

Let us make the change of variable $t \mapsto  T-t$ and, setting 
$x(t):= \overleftarrow{\xi}(T-t)$,    we infer that 

\begin{equation} \label{e:} \left\{ \begin{array}{ll} 
(i) & x(T-\Omega ) \; = \;  s \; \in C(T-\Omega ) \;\mbox{\rm 
and}\; x(T)=x\\ 
(ii) & \forall \; t \in [T-\Omega,T], \; \;  x'(t) \; \in \; F(t,t-(T-\Omega),x(t))\\
(iii) & \forall \; t \in [T-\Omega,T], \; \; x(t) \; \in \; K(t,t-(T-\Omega))\\
\end{array} \right. 
\end{equation}
This means that $(T,x,\Omega ,s)$ belongs to the graph of the 
Cournot map $\mbox{\rm Cour}_{F}(\mathbb{K},\mathbb{C})$.
 \hfill $\;\; \blacksquare$ 
\vspace{ 5 mm} 
 
\subsection{Regulation of Viable Evolutions}
 
Denote by $T^{\star \star}_{K}(x)$ the closed convex hull (or the 
bipolar) of the tangent cone $T_{K}(x)$ to $K$ at $x\in K$.

Recall\footnote{See \emph{Set-valued analysis}, \cite[Aubin \& 
Frankowska]{af90sva}, \emph{Variational Analysis}, 
\cite[Rockafellar \& Wets]{RockafellarWets} and Chapter 18, p. 
713, of \emph{Viability Theory.  New Directions}, \cite[Aubin, 
Bayen \&  Saint-Pierre]{absp}.} that the \index{(forward) 
derivative} \emph{(forward) convexified derivative} $D^{\star 
\star}\mathbb{V}(t,x)$ of a tube $\mathbb{V}$ is defined by

\begin{equation} \label{e:}   
\mbox{\rm Graph}(D^{\star \star}\mathbb{V}(t,x)) \; := \;  
T^{\star \star}_{\mbox{\rm Graph}(\mathbb{V})}(t,x)
\end{equation}
Hence the graph of the convexified forward derivative is a closed 
convex cone (therefore, a set-valued map analogue of a linear 
operator, called a closed convex process in \emph{Convex 
analysis}, \cite[Rockafellar]{r70ca}).

We introduce the concept of regulation map:

\begin{Definition} 
\symbol{91}\textbf{Regulation Map}\symbol{93}\label{} Let us 
consider a set-valued map $F:(t,d,x) \in \mbox{\rm 
Graph}(\mathbb{K})   \leadsto F(t,d,x) \subset X$ and a tube  
$\mathbb{V}: (t,x) \in \mbox{\rm Graph}(\mathbb{K})  \leadsto 
\mathbb{V}(t,x)$. The \index{regulation map} \emph{regulation 
map} $R_{(F,\mathbb{V})}: (t,d,x,s)  \leadsto 
R_{(F,\mathbb{V})}(t,d,x,s)$ is defined by 

\begin{equation} \label{e:}   
R_{(F,\mathbb{V})}(t,d,x,s) \; := \; \left\{ u \in  F(t,d,x)  
\mbox{ such that} \; 0 \; \in \;   D^{\star 
\star}\mathbb{V}(t,x,d,s)(1,u, 1)\right\}  
\end{equation}

\end{Definition}

\textbf{Remark} --- \hspace{ 2 mm} Observe that if $\mathbb{V}$ 
is a single-value differentiable map, then 
$R_{(F,\mathbb{V})}(t,d,x,s)$ is the subset of directions $u \in 
F(t,d,x) $ such that

\begin{displaymath}   
0 \; = \;  \frac{\partial \mathbb{V}(t,d,x,s)}{\partial t} + 
\frac{\partial \mathbb{V}(t,d,x,s)}{\partial d} + \left\langle 
\frac{\partial \mathbb{V}(t,d,x,s)}{\partial x},u  \right\rangle 
\end{displaymath}
which is a McKendrik partial differential equation of 
age-structure problems (see \cite[Aubin]{RegulationBirths-11}). 
\hfill $\;\; \blacksquare$ \vspace{ 5 mm}

Therefore, one can reformulate the Viability Theorem in this 
framework:

\begin{Theorem} 
\symbol{91}\textbf{The Viability 
Theorem}\symbol{93}\label{}\index{} Let us assume that $F$ is 
Marchaud and that the departure and environmental tubes are 
closed. Denote by $\mathbb{V}:= \mbox{\rm 
Cour}_{F}(\mathbb{K},\mathbb{C})$ the Cournot map. Then its 
graph   is closed. The viable evolutions  $x(\cdot) \in 
\mathcal{A}_{F}(T,x)[\Omega,s]$  starting from $s$ at time 
$T-\Omega $ and arriving at $x$ at time $T$ are governed by the 
following differential inclusion involving the regulation map: 

\begin{equation} \label{e:}   
\forall \; t \in [T-\Omega(T,x), T], \; \; x'(t) \; \in \; 
R_{(F,\mathbb{V})}(t,t-(T-\Omega),x(t),s)
\end{equation}
\end{Theorem}

\textbf{Proof} --- \hspace{ 2 mm} The Viability Theorem states 
that whenever the map $F$ is Marchaud,  the capture basin is the 
largest set of elements $(T,\Omega ,x,s )$ between $  
\mathcal{C}$ and $\mathcal{K}$ which is  closed and locally 
viable.  Therefore, the backward velocities $\overleftarrow{u} 
\in F(t,d,x)$   such that $(-1,-1,-\overleftarrow{u},0) \in -( 
\{1\}\times \{1\}\times F(t,d,x)  \times \{0\})$  belongs to 
convexified tangent cone to the viable capture basin $\mbox{\rm 
Capt}_{(\ref{e:LX-DifIncl2})}(\mathcal{K},\mathcal{C})$. By 
Theorem~\ref{t:VbyChaCournot}, p. \pageref{t:VbyChaCournot}, this 
means that  $(-1,-\overleftarrow{u},-1,0) $ belongs to the 
convexified tangent cone $T^{\star \star}_{\mbox{\rm 
Graph}(\mathbb{V})}(\tau,\xi, \omega, \sigma)$ to the graph of 
$\mathbb{V}$. Recalling that $T^{\star \star}_{\mbox{\rm 
Graph}(\mathbb{V})}(\tau,\xi, \omega, \sigma) =\mbox{\rm 
Graph}(D^{\star \star}\mathbb{V}(\tau,\xi, \omega, \sigma) )$, we 
infer that $0 \in - F(\tau, \omega, \xi)\cap D^{\star 
\star}\mathbb{V} (\tau,\xi, \omega, \sigma) (-1, 
-\overleftarrow{u}, -1)$. Therefore, the forward directions 
   $u:=-\overleftarrow{u}$ belongs to 
$R_{(F,\mathbb{V})}(t,d,x,s) $, so that the forward velocities 
$x'(t)$ which regulate the forward viable evolutions $x(\cdot) 
\in \mathcal{A}_{F}(\Omega ,T,x) $ starting at $s$ are the ones 
which belong to $F(t,t-(T-\Omega),x(t))$ and satisfy  $0 \in 
D^{\star \star}\mathbb{V}(t,x(t),t-(T-\Omega),s)(1,x'(t),1)$, 
i.e., which belong to the regulation map 
$R_{(F,\mathbb{V})}(t,t-(T-\Omega),x(t),s)$. 
 \hfill $\;\; \blacksquare$ \vspace{ 5 mm}

\section{Hamilton-Jacobi-Cournot-McKendrik Optimization Problem}


The Cournot map $\mbox{\rm Cour}_{F}(\mathbb{K},\mathbb{C})$ may 
contain \emph{more than one viable evolution}  arriving at $(T,x) 
\in \mbox{\rm Graph}(\mathbb{K})$, starting, for example, from 
the Cournot beginning $T- \Omega 
_{F}(\mathbb{K},\mathbb{C})(T,x)$ and arriving at $x$ at arrival 
time $T$. Hence the question of selecting viable evolutions 
arises.

There are two classes of selection procedures for reducing this 
set of evolutions. The first one is to operate at each time at 
the level of the \index{regulation map} \emph{regulation map} by 
selection one control in  $R_{(F,\mathbb{V})}: (t,d,x,s)  
\leadsto R_{(F,\mathbb{V})}(t,d,x,s)$ (see for instance Section 
11.3.1, p. 453, of \emph{Viability Theory.  New Directions}, 
\cite[Aubin, Bayen \&  Saint-Pierre]{absp}).   

The other class of selection procedures in an intertemporal one, 
which consists in using an intertemporal cost functional on 
evolutions $x(\cdot) \in \mathcal{A}_{F}(T,x)$ (depending for 
instance on departure cost functions and velocity dependent cost 
functions)   and looking for the viable evolutions which minimize 
this intertemporal criterion.  

\begin{Definition} 
\symbol{91}\textbf{Departure Cost Functions and 
Lagrangian}\symbol{93}\label{d:DepArrCost}\index{} We consider   
two ``cost functions''  $\mathbf{c}$ and $\mathbf{l}$: 

\begin{enumerate}
\item a \index{instantaneous condition}\emph{departure cost 
condition} function $ (d,s) \mapsto \mathbf{c}(d,s) \in 
\mathbb{R}\cup \{+\infty\} $;

\item an \index{intertemporal cost functional} 
\emph{intertemporal cost functional}, called in short a 
\emph{Lagrangian}, $\mathbf{l}: (t,d,x,u) \mapsto 
\mathbf{l}(t,d,x,u) \in \mathbb{R}\cup \{+\infty \}$  $u \mapsto 
\mathbf{l}(t,d,x,u)$ depending on time, duration, state and 
velocity;
 \end{enumerate}
with which we associate 
\begin{enumerate}   
\item the departure tube $C: \mathbb{R}^{}  \leadsto X$ defined by

\begin{displaymath}   
C(d) \; := \;  \left\{ s \in X \; \mbox{ such that} \; 
\mathbf{c}(d,s) \; <\; +\infty \right\}
\end{displaymath}   

\item  the set-valued map $F_{\mathbf{l}}:\mathbb{R}^{} \times 
\mathbb{R}^{}_{+} \times   X  \leadsto X$ defined by

\begin{equation} \label{e:}   
F_{\mathbf{l}}(t,d,x) \; := \;  \left\{ u \in X \; \mbox{ such 
that} \; \mathbf{l}(t,d,x,u) \; <\; +\infty \right\}
\end{equation}
and the \index{arrival map} \emph{arrival map} 
$\mathcal{A}_{\mathbf{l}}(T,x)[\Omega,s ]$  associating with any 
final pair $(T,x)$ the set of evolutions $x(\cdot)$ governed by 
the \index{differential inclusion} \emph{differential inclusion} 
\begin{equation} \label{e:LX-DifIncl00}
\forall \;  t \in [T-\Omega ,T], \; \;  x'(t) \; \in \;  
F_{\mathbf{l}}(t,t-(T-\Omega),x(t)) 
\end{equation}
starting at $s\in C(T-\Omega )$ and arriving at the terminal 
condition $x(T) \; = \; x$ at time $T$.
\end{enumerate}
 
\end{Definition}

We have to define the intertemporal cost functional. We begin by 
the simpler case when no constraint function is taken into 
account. 

\begin{Definition} 
\symbol{91}\textbf{The Hamilton-Jacobi-Cournot-McKendrik 
Valuation Function}\symbol{93}\label{}\index{} We associate with 
the data defined in Definition~\ref{d:DepArrCost}, 
p.\pageref{d:DepArrCost} the 
\index{Hamilton-Jacobi-Cournot-McKendrik valuation function} 
\emph{Hamilton-Jacobi-Cournot-McKendrik valuation function} 
$V_{l}(\mathbf{c})$    defined by

\begin{equation} \label{e:} 
V_{l}(\mathbf{c})(T,x) \; := \; \inf_{\Omega \geq 0} 
\inf_{x(\cdot) \in \mathcal{A}_{\mathbf{l}}(T,x)}\; \left( 
\mathbf{c}(T-\Omega ,x(T-\Omega )) + \int_{T-\Omega }^{T} 
\mathbf{l}(t,T-(T-\Omega ), x(t)) dt\right)
\end{equation}
\end{Definition} 

The question arises to know wether the infimum 
$V_{l}(\mathbf{c})(T,x) $ is achieved and what are the standard 
properties of the valuation function, and, in particular, what is 
the Hamilton-Jacobi-McKendrik to which it is a solution.

The way to achieve this program is to observe that the epigraph 
of the valuation function is the Cournot map of an auxiliary 
problem we now define (the vertical arrows symbolize this 
property and the fact that these Cournot maps are related to 
intertemporal minimization problems).

We introduce the 
\begin{enumerate}   

\item   auxiliary target  
\begin{displaymath}   
\mathcal{C}_{\uparrow} \; := \; (t,d,x,y,s) \; \mbox{ such that} 
\; \mathbf{c}(t,s) <+\infty, \; d=0 \;\mbox{\rm and}\; s=x
\end{displaymath} 

\item  the right hand side 
\begin{equation} \label{e:}   
\mathcal{F}_{\uparrow}(t,d,x,y,s) \; = \; \left\{\{1\} \times 
\{1\} \times  {\cal E}p(\mathbf{l}) \times \{0\}  \right\} 
\end{equation}
of the differential inclusion  \begin{equation} \label{e:}   
t'=1, \; d'=1, \;  x'=u, \; y'\leq \mathbf{l}(t,d,x,u), \; s'=0 
\end{equation}

\item the  Cournot map  $(T,x,y)  \leadsto  \mbox{\rm 
Cour}_{\mathcal{F}_{\uparrow}}(\mathbb{K}_{\uparrow}, 
\mathbb{C}_{\uparrow})(T,x,y)$.
\end{enumerate}

\begin{Definition} 
\symbol{91}\textbf{The Viability Solution to the 
Hamilton-Jacobi-Cournot-McKendrik Optimization 
Problem}\symbol{93}\label{}\index{} We consider the extended 
Cournot map 
\begin{equation} \label{e:}   
(T,x,y)  \leadsto  \mbox{\rm Cour}_{\mathcal{F}_{\uparrow}}( 
\mathbb{K}_{\uparrow}, \mathbb{C}_{\uparrow}) (T,x,y)
\end{equation}
associating with elements $(T,x,y)$ the set of apertures $\Omega 
\geq 0$ and initial values $s=x(T-\Omega )$ of evolutions 
$x(\cdot) \in \mathcal{A}_{\mathbf{l}}(T,x)[\Omega,s ]$ starting 
from $s \in \mathcal{C}(T-\Omega )$.

The \index{Hamilton-Jacobi-Cournot-McKendrik viability solution} 
\emph{Hamilton-Jacobi-Cournot-McKendrik viability solution} 
$W_{l}(\mathbf{c})$ is defined by

\begin{equation} \label{e:}   
W_{l}(\mathbf{c})(T,x) \; := \; \inf_{(T,x,y) \in \mbox{\rm 
Dom}(\mbox{\rm 
Cour}_{\mathcal{F}_{\uparrow}}(\mathbb{K}_{\uparrow}, 
\mathbb{C}_{\uparrow})) } y
\end{equation}
\end{Definition}

As expected, these two functions coincide.

\begin{Theorem} 
\symbol{91}\textbf{The Hamilton-Jacobi-Cournot-McKendrik 
Valuation Function and Viability Solution Coincide}\symbol{93} 
\label{}\index{} 

\begin{equation} \label{e:}   
\forall \;  (T,x), \; \;V_{l}(\mathbf{c})(T,x) \; = \; 
W_{l}(\mathbf{c})(T,x)
\end{equation}
Therefore, the valuation function inherits the properties of 
Cournot maps.
\end{Theorem}

\textbf{Proof} --- \hspace{ 2 mm} Let $(T,x,y) \in \mbox{\rm 
Dom}(\mbox{\rm 
Cour}_{\mathcal{F}_{\uparrow}}(\mathbb{K}_{\uparrow}, 
\mathbb{C}_{\uparrow}))$ and $(\Omega, x(T-\Omega)  )$ belong to  
$\mbox{\rm Cour}_{\mathcal{F}_{\uparrow}}(\mathbb{K}_{\uparrow}, 
\mathbb{C}_{\uparrow})(T,x,y)$. Then there exist $s \in 
C(T-\Omega )$ and $x(\cdot) \in 
\mathcal{A}_{\mathbf{l}}(T,x)[\Omega,s ]$ such that  
$s=x(T-\Omega )  $ and $y(T-\Omega ) \geq  \mathbf{c}(T-\Omega, 
x(T-\Omega ) )$. Since

\begin{equation} \label{e:}   
y(T-\Omega ) \; \leq  \;  y - \int_{T-\Omega 
}^{T}\mathbf{l}(t,t-(T-\Omega ),x(t), u(t))dt
\end{equation}    
because $y'(t) \leq  \mathbf{l}(t,t-(T-\Omega ),x(t), u(t))$ and 
$y(T)=y$, we infer that

\begin{equation} \label{e:}   
\mathbf{c}(T-\Omega, x(T-\Omega ) ) + \int_{T-\Omega 
}^{T}\mathbf{l}(t,t-(T-\Omega ),x(t), u(t))dt  \; \leq \; y
\end{equation}

By taking the infimum over $\Omega \geq 0$ and, next, over the  
$x(\cdot) \in \mathcal{A}_{\mathbf{l}}(T,x)[\Omega ]$, we deduce 
that the valuation function $V_{l}(\mathbf{c})(T,x) \leq y$. By 
taking the infimum over the set of  $y$ satisfying $(T,x,y) \in 
\mbox{\rm Dom}(\mbox{\rm 
Cour}_{\mathcal{F}_{\uparrow}}(\mathbb{K}_{\uparrow}, 
\mathbb{C}_{\uparrow}))$, we obtain inequality 
$V_{l}(\mathbf{c})(T,x) \leq W_{l}(\mathbf{c})(T,x)$.

For proving the opposite inequality, let us fix $\varepsilon >0$ 
and choose $\Omega _{\varepsilon} \geq 0$ and an evolution $x 
_{\varepsilon}
 (\cdot) \in \mathcal{A}_{\mathbf{l}}(T,x)[\Omega _{\varepsilon} 
]$ such that 

\begin{equation} \label{e:}   
\mathbf{c}(T-\Omega_{\varepsilon} 
,x_{\varepsilon}(T-\Omega_{\varepsilon} )) + 
\int_{T-\Omega_{\varepsilon} }^{T} 
\mathbf{l}(t,T-(T-\Omega_{\varepsilon} ), x_{\varepsilon}(t)) dt  
\; \leq \;  V_{l}(\mathbf{c})(T,x) +\varepsilon
\end{equation}
Let us set

\begin{equation} \label{e:}   
y_{\varepsilon}(t) \; := \; V_{l}(\mathbf{c})(T,x) +\varepsilon 
-\int_{t }^{T} \mathbf{l}(\tau,\tau-(T-\Omega_{\varepsilon} ), 
x_{\varepsilon}(\tau), u_{\varepsilon}(\tau)) dt 
\end{equation}

We thus observe that 
$(x_{\varepsilon}(\cdot),y_{\varepsilon}(\cdot))$ is a solution 
to the differential inclusion 
$(x_{\varepsilon}'(\cdot),y'_{\varepsilon}(\cdot)) \; \in \; 
{\cal E}p(\mathbf{l})$,     that $x_{\varepsilon}(T)=x$ and that 
$y_{\varepsilon}(T)= V_{l}(\mathbf{c})(T,x) +\varepsilon $ and 
$y(T-\Omega )=V_{l}(\mathbf{c})(T,x) +\varepsilon-\int_{T-\Omega  
}^{T} \mathbf{l}(\tau,\tau-(T-\Omega_{\varepsilon} ), 
x_{\varepsilon}(\tau), u_{\varepsilon}(\tau)) dt \; \leq \;  
V_{l}(\mathbf{c})(T,x) +\varepsilon$. Hence $(T,x, 
V_{l}(\mathbf{c})(T,x) +\varepsilon)$ belongs to the domain of 
the auxiliary Cournot map. This implies that 
$W_{l}(\mathbf{c})(T,x) \leq  V_{l}(\mathbf{c})(T,x) 
+\varepsilon$. Letting $\varepsilon  \rightarrow  0+$ implies 
$W_{l}(\mathbf{c})(T,x) \leq  V_{l}(\mathbf{c})(T,x) $ and thus 
the equality we were looking for. \hfill $\;\; \blacksquare$ 
\vspace{ 5 mm} 


We  recall that the viability solution, when it is 
differentiable, is a solution to the Hamilton-Jacobi equation 
satisfying the trajectory conditions. Otherwise, when it is not 
differentiable, but only lower semicontinuous, we can give a 
meaning to a solution as a solution  in the 
Barron-Jensen/Frankoska sense, using for that purpose 
subdifferential of lower semicontinuous functions defined in 
non-smooth analysis (\emph{Set-valued analysis}, \cite[Aubin \& 
Frankowska]{af90sva}, \cite[Aubin, Bayen, 
Saint-Pierre]{absp06hj,absp}). This is not that important for two 
reasons: all other properties of viability solutions that are 
proven in this paper are derived directly from the properties of 
capture basins without using the concept of derivatives, usual or 
generalized. In particular, the fact that the viability solution 
is a solution to the Hamilton-Jacobi equation derives from the 
tangential conditions characterization of viable-capture basins 
provided by the Viability Theorem. 

The adaptation of the results of Chapters~13 and 17 of 
\emph{Viability Theory.  New Directions}, \cite[Aubin, Bayen \&  
Saint-Pierre]{absp}  is straightforward.

\clearpage \tableofcontents 
\end{document}